\documentclass[11pt]{amsart}
\usepackage{amssymb}
\usepackage{amsmath}
\usepackage{fancyhdr}

\title[Decompositions of complete uniform hypergraphs into Berge cycles]{Decompositions of complete uniform hypergraphs into Hamilton Berge cycles}

\date{\today}
\author{Daniela K\"uhn and Deryk Osthus}
\thanks{The research leading to these results was partially supported by the European Research Council
under the European Union's Seventh Framework Programme (FP/2007--2013) / ERC Grant
Agreements no. 258345 (D.~K\"uhn) and 306349 (D.~Osthus).} 

\newtheorem{firstthm}{Proposition}
\newtheorem{theorem}[firstthm]{Theorem}

\newtheorem{lemma}[firstthm]{Lemma}

\newtheorem{conj}[firstthm]{Conjecture}

\def\noproof{{\unskip\nobreak\hfill\penalty50\hskip2em\hbox{}\nobreak\hfill%
       $\square$\parfillskip=0pt\finalhyphendemerits=0\par}\goodbreak}
\def\endproof{\noproof\bigskip}
\newdimen\margin   % needed for macros \textdisplay & \ltextdisplay
\def\textno#1&#2\par{%
   \margin=\hsize
   \advance\margin by -4\parindent
          \setbox1=\hbox{\sl#1}%
   \ifdim\wd1 < \margin
      $$\box1\eqno#2$$%
   \else
      \bigbreak
      \hbox to \hsize{\indent$\vcenter{\advance\hsize by -3\parindent
      \it\noindent#1}\hfil#2$}%
      \bigbreak
   \fi}
\def\proof{\removelastskip\penalty55\medskip\noindent{\bf Proof. }}

\begin{document}

\def\COMMENT#1{}
\def\TASK#1{}

\def\eps{{\varepsilon}}
\newcommand{\ex}{\mathbb{E}}
\newcommand{\pr}{\mathbb{P}}
\newcommand{\cB}{\mathcal{B}}
\newcommand{\cS}{\mathcal{S}}
\newcommand{\cF}{\mathcal{F}}
\newcommand{\cC}{\mathcal{C}}
\newcommand{\cP}{\mathcal{P}}
\newcommand{\cQ}{\mathcal{Q}}
\newcommand{\cR}{\mathcal{R}}
\newcommand{\cK}{\mathcal{K}}
\newcommand{\cD}{\mathcal{D}}
\newcommand{\cI}{\mathcal{I}}
\newcommand{\cV}{\mathcal{V}}
\newcommand{\1}{{\bf 1}_{n\not\equiv \delta}}
\newcommand{\eul}{{\rm e}}

\begin{abstract}  \noindent
In 1973 Bermond, Germa, Heydemann and Sotteau conjectured that if $n$ divides $\binom{n}{k}$, 
then the complete $k$-uniform hypergraph on $n$ vertices  has a decomposition into Hamilton Berge cycles.
Here a Berge cycle consists of an alternating sequence $v_1,e_1,v_2,\dots,v_n,e_n$
of distinct vertices $v_i$ and distinct edges $e_i$ so that each $e_i$ contains $v_i$ and $v_{i+1}$.
So the divisibility condition is clearly necessary.
In this note, we prove that the conjecture holds whenever $k \ge 4$ and $n \ge 30$.
Our argument is based on the Kruskal-Katona theorem.
The case when $k=3$ was already solved by Verrall, building on results of Bermond.
\end{abstract}
\maketitle

\section{Introduction}\label{intro}
A classical result of Walecki~\cite{lucas} states that the complete graph $K_n$ on $n$ vertices has
a Hamilton decomposition if and only if $n$ is odd.
(A Hamilton decomposition of a graph $G$ is a set of edge-disjoint Hamilton cycles containing all edges of $G$.)
Analogues of this result were proved for complete digraphs by Tillson~\cite{till} and more recently for (large) tournaments in~\cite{Kelly}.
Clearly, it is also natural to ask for a hypergraph generalisation of Walecki's theorem.

There are several notions of a hypergraph cycle, the earliest one is due Berge:
A \emph{Berge cycle} consists of an alternating sequence $v_1,e_1,v_2,\dots,v_n,e_n$
of distinct vertices $v_i$ and distinct edges $e_i$ so that each $e_i$ contains $v_i$ and $v_{i+1}$.
A Berge cycle is a Hamilton (Berge) cycle of a hypergraph $G$ if $\{v_1,\dots,v_n\}$ is the vertex set of~$G$
and each $e_i$ is an edge of~$G$.
So a Hamilton Berge cycle has $n$ edges.

Let $K_n^{(k)}$ denote the complete $k$-uniform hypergraph on $n$ vertices. 
Clearly, a necessary condition for the existence of 
a decomposition of $K_n^{(k)}$ into Hamilton Berge cycles is that $n$ divides $\binom{n}{k}$.
Bermond, Germa, Heydemann and Sotteau~\cite{BGHS} conjectured that this condition is also sufficient.
For $k=3$, this conjecture follows by combining the results of Bermond~\cite{bermond} and Verrall~\cite{verrall}. 

We show that as long as $n$ is not too small, the conjecture holds for $k\ge 4$ as well.
\begin{theorem} \label{simple}
Suppose that $4 \le k <n$, that $n \ge 30$ and that $n$ divides $\binom{n}{k}$.
Then the complete $k$-uniform hypergraph 
$K_n^{(k)}$ on $n$ vertices  has a decomposition into Hamilton Berge cycles.
\end{theorem}
Recently, Petecki~\cite{petecki} considered a restricted type of decomposition into Hamilton Berge cycles and determined those $n$ for which $K_n^{(k)}$
has such a restricted decomposition.

Walecki's theorem has a natural extension to the case when $n$ is even:
in this case, one can show that $K_n-M$ has a Hamilton decomposition, 
whenever $M$ is a perfect matching. 
Similarly, the results of Bermond~\cite{bermond} and Verrall~\cite{verrall} together imply that
for all $n$, either $K_n^{(3)}$ or $K_n^{(3)}-M$ have a decomposition into Hamilton Berge cycles.

We prove an analogue of this for $k \ge 4$. Note that Theorem~\ref{thm:Bergedec} immediately implies 
Theorem~\ref{simple}.

\begin{theorem}\label{thm:Bergedec}
Let $k,n\in \mathbb{N}$ be such that $3\le k< n$.
\begin{itemize}
\item[{\rm (i)}] Suppose that $k\ge 5$ and $n\ge 20$ or that $k=4$ and $n\ge 30$. Let $M$ be any set consisting of less than $n$ edges of  $K^{(k)}_n$ such that
$n$ divides $|E(K^{(k)}_n)\setminus M|$. Then $K^{(k)}_n-M$ has a decomposition into Hamilton Berge cycles.
\item[{\rm (ii)}] Suppose that $k=3$ and $n\ge 100$. If $\binom{n}{3}$ is not divisible by $n$, let $M$ be any perfect matching%
    \COMMENT{If $\binom{n}{3}$ is not divisible by $n$ then $3$ divides $n$ (since $\binom{n}{3}=\frac{n}{3}\binom{n-1}{2}$).
So $K^{(3)}_n$ has a perfect matching. But also, in this case $n$ will divide $\binom{n}{3}-\frac{n}{3}=\frac{n}{3}(\binom{n-1}{2}-1)$
since $\binom{n-1}{2}\equiv 1\mod 3$)}
in $K^{(k)}_n$,
otherwise let $M:=\emptyset$. Then $K^{(3)}_n-M$ has a decomposition into Hamilton Berge cycles.
\end{itemize}
\end{theorem}

Note that if $k$ is a prime and $\binom{n}{k}$ is not divisible by $n$, then $k$ divides $n$ and so in this case one can take the set $M$ in (i)
to be a union of perfect matchings.%
\COMMENT{this is not generally possible: e.g. if $n=22$ and $k=4$}
Also note that (ii) follows from the results of~\cite{bermond,verrall}.
However, our proof is far simpler, so we also include it in our argument.

%There are several other natural notions of hypergraph Hamilton cycles.
%Ber\-mond, Germa, Heydemann and Sotteau~\cite{BGHS} defined a \emph{Berge $\ell$-cycle} to be a Berge cycle 
%where consecutive edges intersect in exactly $\ell$ vertices.
%They also made the following conjecture.
%\begin{conj} \label{typeconj} Let $k,\ell,n\in \mathbb{N}$ be such that $3\le k< n$.
%Suppose that $n$ divides $\binom{n}{k}$ and that $1 \le \ell <k$.
%Then the complete $k$-uniform hypergraph 
%$K_n^{(k)}$ on $n$ vertices  has a decomposition into Hamilton Berge $\ell$-cycles.
%\end{conj}
%This clearly implies Theorem~\ref{simple} and is wide open.
Another popular notion of a hypergraph cycle is the following: 
a $k$-uniform hypergraph~$C$ is an $\ell$-\textit{cycle} if there 
exists a cyclic ordering of the vertices of~$C$ such that every edge of~$C$ 
consists of~$k$ consecutive vertices and such that every pair of {consecutive} edges 
(in the natural ordering of the edges) intersects in precisely $\ell$ vertices. 
If $\ell=k-1$, then $C$ is called a \emph{tight cycle} and if $\ell=1$, then $C$ is called a \emph{loose cycle}.
We conjecture an analogue of Theorem~\ref{simple} for Hamilton $\ell$-cycles.

\begin{conj} \label{lcycles}
For all $k,\ell\in \mathbb{N}$ with $\ell<k$ there exists an integer $n_0$ such that the following
holds for all $n\ge n_0$. Suppose that $k-\ell$ divides $n$ and that $n/(k-\ell)$ divides $\binom{n}{k}$.
Then $K_n^{(k)}$ has a decomposition into Hamilton $\ell$-cycles.
\end{conj}
To see that the divisibility conditions are necessary, note that every Hamilton $\ell$-cycle contains exactly $n/(k-\ell)$ edges.
%There are several necessary conditions which need to be satisfied for Conjecture~\ref{lcycles} to hold:
%Firstly, every Hamilton $\ell$-cycle contains exactly $n/(k-\ell)$ edges, so we need $n/(k-\ell) \in \mathbb{N}$.
Moreover, it is also worth noting the following: consider the number $N:=\frac{k-\ell}{n}\binom{n}{k}$ of cycles we require in the decomposition.
%Then we need $N \in \mathbb{N}$ (as it corresponds to the number of 
%cycles in the decomposition). 
The divisibility conditions ensure that $N$ is not only an integer but also 
a multiple of $f:=(k-\ell)/h$, where $h$ is the highest common factor of $k$ and $\ell$.%
    \COMMENT{Indeed, $N=\frac{k-\ell}{k}\binom{n-1}{k-1}$. Let $A:=kN=(k-\ell)\binom{n-1}{k-1}$. Then both $k$ and $k-\ell$ divide $A$.
Thus $k(k-\ell)/hcf(k-\ell,k)$ divides $A$. But $hcf(k-\ell,k)=hcf(\ell,k)=h$. So $N/f=Ah/(k(k-\ell))\in\mathbb{N}$.}
This is relevant as one can construct a regular hypergraph from the edge-disjoint union of $t$ edge-disjoint 
Hamilton $\ell$-cycles if and only if $t$ is a multiple of $f$.%
    \COMMENT{Write $k=a(k-\ell)+b$ with $0 \le b<k-\ell$. 
Denote the edges of a cycle $C$ by $e_i$.
 Consider the final segment $S$ of $k-\ell$ vertices of an edge $e_1$ on a cycle. (Note its degree sequence is identical to that of the $k-\ell$-segment preceding $S$  etc.).
Note $S$ is completely contained in $e_1\dots,e_a$.
Let $S_b$ denote the final $b$ vertices of $S$. Then  $S \cap e_{a+1}=S_b$ (and no other edges $e_i$ of the cycle meet $S$).
By def of $f$, we can split $S$ into $f$ successive subsegments of size $(k-\ell)/f=h$. 
Now note $hcf(k-\ell,b)=hcf(k-\ell,k-a(k-\ell))=hcf(k-\ell,k)=hcf(\ell,k)=h$.
So $b$ is also a multiple of $h$ and the above also splits $S_b$ into successive subsegments of size $h$.
Now take $f$ `shifts' of $C$ (including the identity), where each shift starts $h$ vertices later than the previous one.
Then together these cycles form a regular hypergraph.
Note this must actually be a subhypergraph of the tight cycle.}
%To see that the divisibility conditions are necessary, note 

The `tight' case $\ell=k-1$ of Conjecture~\ref{lcycles} was already formulated by Bailey and Stevens~\cite{baileystevens}.
In fact, if $n$ and $k$ are coprime, the case $\ell=k-1$ already corresponds to a conjecture made independently by Baranyai~\cite{bconj} and Katona on so-called `wreath decompositions'.
A $k$-partite analogue of the `tight' case of Conjecture~\ref{lcycles} was recently proved by Schroeder~\cite{schroeder}.

Conjecture~\ref{lcycles} is known to hold `approximately' (with some additional additional divisibility conditions on $n$), i.e.~one can find a set of edge-disjoint Hamilton $\ell$-cycles which together
cover almost all the edges of $K_n^{(k)}$. This is a very special case of results in~\cite{BF,FK,FKL} which guarantee approximate decompositions of quasi-random
uniform hypergraphs into Hamilton $\ell$-cycles (again, the proofs need $n$ to satisfy additional divisibility constraints).

\section{Proof of Theorem~\ref{thm:Bergedec}}

Before we can prove Theorem~\ref{thm:Bergedec} we need to introduce some notation. Given integers
$0\le k\le n$, we will write $[n]^{(k)}$ for the set consisting of all
$k$-element subsets of $[n]:=\{1,\dots,n\}$. 
The \emph{colexicographic order} on $[n]^{(k)}$ is the order in which $A<B$ if and only if the largest element of $(A\cup B)\setminus (A\cap B)$ lies in~$B$
(for all distinct $A,B\in [n]^{(k)}$). The \emph{lexicographic order}
on $[n]^{(k)}$ is the order in which $A<B$ if and only if the smallest element of
$(A\cup B)\setminus (A\cap B)$ lies in~$A$.
Given $\ell\in\mathbb{N}$ with $\ell\le k$ and a set $S\subseteq [n]^{(k)}$, the \emph{$\ell$th lower shadow of $S$} is the set
$\partial^-_\ell(S)$ consisting of all those $t\in [n]^{(k-\ell)}$ for which there exists $s\in S$ with $t\subseteq s$.
Similarly, given $\ell\in\mathbb{N}$ with $k+\ell\le n$ and a set $S\subseteq [n]^{(k)}$, the \emph{$\ell$th upper shadow of $S$} is the set
$\partial^+_\ell(S)$ consisting of all those $t\in [n]^{(k+\ell)}$ for which there exists $s\in S$ with $s\subseteq t$.
We need the following consequence of the Kruskal-Katona theorem~\cite{Katona,Kruskal}.

\begin{lemma}\label{thm:KK} $\ $
\begin{itemize}
\item[{\rm (i)}] Let $k,n\in \mathbb{N}$ be such that $3\le k\le n$. Given a nonempty $S\subseteq [n]^{(k)}$, define $s\in\mathbb{R}$ by $|S|=\binom{s}{k}$. Then $|\partial^-_{k-2}(S)|\ge \binom{s}{2}$.
\item[{\rm (ii)}] Suppose that $S'\subsetneq [n]^{(2)}$ and let $c,d\in \mathbb{N}\cup \{0\}$ be such that
$c<n$, $d< n-(c+1)$ and $|S'|=cn-\binom{c+1}{2}+d$. If $n\ge 100$ and $c\le 8$ then $|\partial^+_1(S')|\ge c\binom{n-c}{2}+2dn/5$.
\item[{\rm (iii)}] If $S'\subseteq [n]^{(2)}$ and $|S'|\le n-1$ then $|\partial^+_{2}(S')|\ge |S'|\binom{n-|S'|-1}{2}+\binom{|S'|}{2}(n-|S'|-1)$.%
    \COMMENT{Could use $\partial^+_{2}(S')\ge |S'|\binom{n-|S'|-1}{2}+\binom{|S'|}{2}(n-|S'|-1)+\binom{|S'|}{3}$ instead. This would probably improve our bound on $n$.}
\end{itemize}
\end{lemma}
\proof
The Kruskal-Katona theorem states that the size of the lower shadow of a set $S\subseteq [n]^{(k)}$ is minimized if $S$ is an initial segment of $[n]^{(k)}$ in the colexicographic order.
(i) is a special case of a weaker (quantitative) version of this due to Lov\'asz~\cite{Lovasz}.
In order to prove (ii) and (iii), note that whenever $A,B\in [n]^{(k)}$ then $A<B$ in the colexicographic order if and only if
$[n]\setminus A<[n]\setminus B$ in the lexicographic order on $[n]^{(n-k)}$ with the order of the ground set reversed.
Thus, by considering complements, it follows from the Kruskal-Katona theorem that the size of the upper shadow of a set $S'\subseteq [n]^{(k)}$ is minimized if $S'$ is an initial
segment of $[n]^{(k)}$ in the lexicographic order. This immediately implies~(iii).
Moreover, if $S'$, $c$ and $d$ are as in~(ii), then%
    \COMMENT{The initial segment of size $|S'|$ in lex is $L:=\{12,\dots,1n,23,\dots,2n,\dots,c(c+1),\dots, cn, (c+1)(c+2),\dots,(c+1)(c+d+1)\}$.
$\partial^+_1 L$ consists of all $3$-sets containing 1 ($\binom{n-1}{2}$ choices), all $3$-sets containing 2 but not 1 ($\binom{n-2}{2}$ choices),
and so on. The term $\binom{n-c}{2}$ counts all 3-sets containing $c$ but none of $1,\dots,c-1$, the term $d(n-c-2)-\binom{d}{2}$ counts all
3-sets containing one of $(c+1)(c+2),\dots,(c+1)(c+d+1)$ but none of $1,\dots,c$ (there are $d$ choices to choose one of $(c+1)(c+2),\dots,(c+1)(c+d+1)$,
and $n-c-2$ choices for the remaining number, however, this counts the 3-sets which consist of $c+1$ and two numbers from $c+2,\dots,c+d+1$ twice)}
\begin{align*}
|\partial^+_1 S'| & \ge \binom{n-1}{2}+\binom{n-2}{2}+\dots + \binom{n-c}{2}+d(n-c-2)-\binom{d}{2}\\
& \ge c\binom{n-c}{2}+\frac{2}{5} dn,
\end{align*}
as%
   \COMMENT{To see the final inequality note that $d(n-c-2)-\binom{d}{2}=d(n-c-2-(d-1)/2)\ge d(n/2-c-2)\ge d(n/2-10)\ge 2dn/5$ since $d\le n$, $c\le 8$ and $n\ge 100$.}
required.
\endproof

We will also use the following result of Tillson~\cite{till} on Hamilton decompositions of complete digraphs.
(The \emph{complete digraph $DK_n$} on $n$ vertices has a directed edge $xy$ between every ordered pair $x\neq y$ of vertices. So $|E(DK_n)|=n(n-1)$.) 

\begin{theorem}\label{thm:tillson}
The complete digraph $DK_n$ on $n$ vertices has a Hamilton decomposition if and only if $n\neq 4,6$.
\end{theorem}

\removelastskip\penalty55\medskip\noindent{\bf Proof of Theorem~\ref{thm:Bergedec}. }
The first part of the proof for (i) and (ii) is identical. So let $M$ be as in (i),(ii).
(For~(ii) note that if $\binom{n}{3}$ is not divisible by $n$, then $3$ divides $n$ and $n$ divides $\binom{n}{3}-\frac{n}{3}$.) Let
$$\ell:=\left\lfloor\frac{\binom{n}{k}-|M|}{n(n-1)}\right\rfloor \ \ \ \ \text{and} \ \ \ \  m:=\frac{\binom{n}{k}-|M|-\ell n(n-1)}{n}.$$
Note that $m<n-1$ and $m\in \mathbb{N}\cup \{0\}$ since $n$ divides $\binom{n}{k}-|M|$.
Define an auxiliary (balanced) bipartite graph $G$ with vertex classes $A_*$ and $B$ of size $\binom{n}{k}-|M|$ as follows. 
Let $A:=E(K^{(k)}_n)$ and $A_*:=A\setminus M$. Let $D_1,\dots,D_\ell$ be copies of the complete digraph $DK_n$ on $n$ vertices.
For each $i\in [\ell]$ let $B_i,B'_i$ be a partition of $E(D_i)$ such that for every pair $xy, yx$ of opposite directed edges, $B_i$
contains precisely one of $xy, yx$.
Apply Theorem~\ref{thm:tillson} to find $m$ edge-disjoint Hamilton cycles $H_1,\dots,H_m$ in
$DK_n$. We view the sets $B_1,\dots,B_\ell$, $B'_1,\dots,B'_\ell$ and $E(H_1),\dots, E(H_m)$ as being pairwise disjoint and
let $B$ denote the union of these sets. So $|B|=|A_*|$. Our auxiliary bipartite graph $G$ contains an edge between $z\in A_*$ and
$xy\in B$ if and only if $\{x,y\}\subseteq z$.

We claim that $G$ contains a perfect matching $F$. Before we prove this claim, let us show how it implies 
Theorem~\ref{thm:Bergedec}. For each $i\in [\ell]$, apply Theorem~\ref{thm:tillson} to obtain a Hamilton decomposition $H^1_i,\dots,H^{n-1}_i$
of $D_i$. For each $i\in [\ell]$ and each $j\in [n-1]$ let $A^j_i\subseteq A$ be the neighbourhood of $E(H^j_i)$ in $F$.
Note that each $A^j_i$ is the edge set of a Hamilton Berge cycle of $K^{(k)}_n-M$.
Similarly, for each $i'\in [m]$ the neighbourhood $A_{i'}$ of $E(H_{i'})$ in $F$ is the edge set of a Hamilton Berge cycle of $K^{(k)}_n-M$.
Since all the sets $A^j_i$ and $A_{i'}$ are pairwise disjoint, this gives a decomposition of $K^{(k)}_n-M$ into Hamilton Berge cycles.

Thus it remains to show that $G$ satisfies Hall's condition. So consider any nonempty set $S\subseteq A_*$ and define $s,a\in \mathbb{R}$ with $k\le s\le n$
and $0<a\le 1$ by $|S|=a\binom{n}{k}=\binom{s}{k}$. Define $b$ by $|N_G(S)\cap B_1|=b\binom{n}{2}$. Note that
$|N_G(S)\cap B_1|\ge \binom{s}{2}$ by Lemma~\ref{thm:KK}(i).
But%
   \COMMENT{Indeed $\frac{\binom{s}{2}^k}{\binom{n}{2}^k}=\left(\frac{s(s-1)}{n(n-1)}\right)^k\ge \frac{((s)_k)^2}{((n)_k)^2}=\frac{ \binom{s}{k}^2}{\binom{n}{k}^2}.$}
$$\frac{b^k}{a^2}\ge \frac{\binom{s}{2}^k\binom{n}{k}^2}{\binom{n}{2}^k\binom{s}{k}^2}\ge 1,$$
and so $b\ge a^{2/k}$. Thus
\begin{align*}
|N_G(S)| &\ge 2\ell |N_G(S)\cap B_1|\ge 2\ell a^{2/k}\binom{n}{2} =a^{2/k}(|B|-|E(H_1)\cup\dots\cup E(H_m)|)\\
& \ge a^{2/k}\left(|A_*|-n(n-2)\right).
\end{align*}
Let%
    \COMMENT{$g<1$ iff $k<n-2$. But if $k\ge n-2$ we don't use (\ref{eq:a}).}
$$g:=\frac{\binom{n}{k}-|A_*|+n(n-2)}{\binom{n}{k}}.$$
So if
\begin{equation}\label{eq:a}
a^{1-2/k}\le \frac{|A_*|}{\binom{n}{k}}-\frac{n(n-2)}{\binom{n}{k}}=1-g,
\end{equation}
then $|N_G(S)|\ge |S|$. We now distinguish three cases.

\medskip

\noindent\textbf{Case~1.} $4 \le k\le n-3$

\smallskip

\noindent Since
$$|A_*|-2n(n-1)\le |A_*|-\left(\binom{n}{k}-|A_*|\right)-2n(n-2)= (1-2g)\binom{n}{k}\le (1-g)^2\binom{n}{k},$$
in this case (\ref{eq:a}) implies that $|N_G(S)|\ge |S|$ if
$|S|\le |A_*|-2n(n-1)$. So suppose that $|S|> |A_*|-2n(n-1)$. Note that if $k\ge 5$ then every $b\in B$ satisfies
$|N_G(b)|=\binom{n-2}{k-2}-|M|\ge \binom{n-2}{3}-n\ge 2n(n-1)$ since $n\ge k+3$ and $n\ge 20$.%
   \COMMENT{Note that $\binom{n-2}{k-2}=\binom{n-2}{3}\frac{n-5}{k-2}\frac{n-6}{k-3}\dots\frac{n-k+1}{4}$. But
$\frac{n-5}{k-2}\frac{n-6}{k-3}\dots\frac{n-k+1}{4}\ge 1$ since $n\ge k+3$.
Also $(n-2)(n-3)(n-4)/6 \ge \frac{16}{6}n^2 \frac{n-2}{n} \frac{n-3}{n} \ge \frac{8}{3} n^2 (1-5/n) \ge  \frac{8}{3} n^2 \frac{3}{4}=2n^2$}
Hence $N_G(S)=B$.

So we may assume that $k=4$ and $S':=B\setminus N_G(S)\neq \emptyset$. Thus $S'_1:=S'\cap B_1\neq \emptyset$
and $|S'|\le (2\ell+2) |S'_1|$. Note that $|N_G(S'_1)|\le |A_*\setminus S|< 2n(n-1)$.
First suppose $|S'_1|\ge 7$. Then $|N_G(S'_1)|\ge 7\binom{n-8}{2}+21(n-8)-|M| > 2n(n-1)$
 by Lemma~\ref{thm:KK}(iii) and our assumption that $n\ge 30$.%
 \COMMENT{$7\binom{n-8}{2}+21(n-8)-n=7/2(n-8)(n-9)+20n-168=7/2(n^2-17n+72)+20n-168
\ge 2n^2+1/2(3 \cdot 30n-7\cdot 17n+7 \cdot 72)+20n -168
=2n^2+n/2(90-70-49+40)+7 \cdot 36 -168 \ge 2n^2$}
So we may assume that $|S'_1|\le 6$. Apply Lemma~\ref{thm:KK}(iii) again to see that
\begin{align*}
|N_G(S')| & \ge |S'_1|\binom{n-7}{2}-|M|\ge \frac{\binom{n-7}{2}}{2\ell+2}|S'|-n\ge \frac{6(n-7)(n-8)}{(n-2)(n-3)+24}|S'|-n \\
& \ge 2|S'|-n>|S'|.
\end{align*}
(Here we use that $|S'|\ge 2\ell > n$%
\COMMENT{$2\ell \ge 2(n-2)(n-3)/24 -4\ge 28(n-3)/12-4 \ge n-28 \cdot 3/12+16 \cdot 27/12-4 > n$}
 and $n\ge 30$.)%
\COMMENT{$6(n-7)(n-8)=6(n^2-15n+56)\ge 6n^2-90 n+300 \ge 2n^2+4 \cdot 30 n-90 n+300 \ge 2(n^2+150) \ge 2(n^2-5n+30)=2((n-2)(n-3)+24)$}
 Thus $|N_G(S)|\ge |S|$, as required.

\medskip

\noindent\textbf{Case~2.} $k= 3$

\smallskip
\noindent Since
$$|A_*|-3n(n-1)\le |A_*|-2\left(\binom{n}{k}-|A_*|\right)-3n(n-2)= (1-3g)\binom{n}{k}\le (1-g)^3\binom{n}{k},$$
in this case (\ref{eq:a}) implies that $|N_G(S)|\ge |S|$ if
$|S|\le |A_*|-3n(n-1)$. So suppose that $|S|> |A_*|-3n(n-1)$
and that $S':=B\setminus N_G(S)\neq \emptyset$. Thus $S'_1:=S'\cap B_1\neq \emptyset$
and $|S'|\le (2\ell+2) |S'_1|\le ((n-2)/3+2)|S'_1|$. Let $c,d\in \mathbb{N}\cup \{0\}$ be such that
$c<n$, $d< n-(c+1)$ and $|S'_1|=cn-\binom{c+1}{2}+d$. Note that $|N_G(S'_1)|\le |A_*\setminus S|< 3n(n-1)$.
Thus $c<8$ since otherwise
$$|N_G(S'_1)|\ge 8\binom{n-8}{2}-|M|\ge 8\binom{n-8}{2}-\frac{n}{3}>\frac{32}{5}\binom{n}{2} > 3n(n-1)$$
by Lemma~\ref{thm:KK}(ii)
and our assumption that $n\ge 100$.%
   \COMMENT{We use that $\binom{n-8}{2}\ge \binom{n}{2}(1-8/n)(1-8/(n-1))\ge  \binom{n}{2}(1-16/(n-1))$ and so
$8\binom{n-8}{2}-n/3=8[\binom{n-8}{2}-n/24]\ge 8\binom{n}{2}(1-\frac{16}{n-1}-\frac{1}{12(n-1)})\ge 8 \binom{n}{2}(1-\frac{17}{n-1})\ge 8\frac{4}{5}\binom{n}{2}$
as $n\ge 100$}
Let $M(S'_1)$ denote the set of all those edges $e\in M$ for which there is a pair $xy\in S'_1$ with $\{x,y\}\subseteq e$.
Thus $M(S'_1)=\partial^+_1(S'_1)\cap M$.
Recall that $M$ is a matching in the case when $k=3$. Thus $|M(S'_1)|\le |S'_1|$. In particular $|M(S'_1)|\le d$ if $c=0$. 
Apply Lemma~\ref{thm:KK}(ii) again to see that%
   \COMMENT{To see the 3rd line, note that if $c\ge 1$ then $\frac{4c}{5}\binom{n}{2}-\frac{n}{3}\ge \frac{4c}{5}\binom{n}{2}-\frac{2cn}{5}= \frac{12c}{10}\frac{n(n-2)}{3}.$
Moreover, if $c=0$ then $\frac{2dn}{5}-d\ge \frac{11d}{10}\frac{n}{3}$.}
\begin{align*}
|N_G(S')| & \ge |N_G(S'_1)| \ge c\binom{n-c}{2}+\frac{2}{5}dn- |M(S'_1)|\\
& \ge \frac{4c}{5}\binom{n}{2}+\frac{2}{5}dn-\begin{cases} n/3 &\mbox{if } c\ge 1 \\
d & \mbox{if } c=0 \end{cases}\\
& \ge (cn+d)\cdot \frac{11}{10}\cdot \frac{n-2}{3}\ge |S'_1|\left(\frac{n-2}{3}+2\right)\ge |S'|,
\end{align*}
where we use that $n\ge 100$. Thus $|N_G(S)|\ge |S|$, as required.

\medskip

\noindent\textbf{Case~3.} $n-2\le k\le n-1$

\smallskip

\noindent If $k=n-1$ then $K_n^{(k)}$ itself is a Hamilton Berge cycle, so there is nothing to show. So suppose that
$k=n-2$. 
In this case, it helps to be more careful with the choice of the Hamilton cycles $H_1,\dots,H_m$:
instead of applying Theorem~\ref{thm:tillson} to find $m$ edge-disjoint Hamilton cycles $H_1,\dots,H_m$ in
$DK_n$, we proceed slightly differently.
Note first that $\ell=0$. Suppose that $n$ is odd. Then $M=\emptyset$ and $m=(n-1)/2$.
If $n$ is even, then $|M|=n/2$ and $m=n/2-1$.
In both cases we can choose $H_1,\dots,H_m$ to be $m$ edge-disjoint Hamilton cycles of~$K_n$.
Then a perfect matching in our auxiliary graph $G$ still corresponds to a decomposition of $K_n^{(k)} -M$ into Hamilton Berge cycles.
%Note that $\ell=0$ and $m$m \ge \frac{1}{n}\binom{n}{2}-1 \ge \frac{n}{3}$.
Also, in both cases $E(H_1)\cup\dots\cup E(H_m)$ contains all but at most $n/2$ distinct elements of $[n]^{(2)}$.

Consider any $b \in B$.
Then%
\COMMENT{
$\ge \binom{n}{2} \left(\frac{n-2}{n}\frac{n-3}{n-1} - \frac{1}{n-1} \right) \ge$}
$$
|N_G(b)| \ge \binom{n-2}{k-2}-|M| = \binom{n-2}{2}-|M| \ge \binom{n}{2}\left( 1-\frac{5}{n-1} \right) \ge \frac{2}{3} \binom{n}{2} \ge \frac{2}{3}|A_*|.
$$
Now consider any $a \in A_*$. 
Then
\COMMENT{$ \binom{k}{2}=\binom{n-2}{2} \ge \binom{n}{2}(1-2\cdot 2/(n-1)) \ge \frac{15}{19}\binom{n}{2}$}
$$
|N_G(a)| \ge \binom{k}{2}-\frac{n}{2} \ge \frac{2}{3} \binom{n}{2} \ge \frac{2}{3}|B|.
$$
So Hall's condition is satisfied and so $G$ has a perfect matching, as required.
\endproof

The lower bounds on $n$ have been chosen so as to streamline the calculations, and could be improved by more careful calculations.

\medskip

{\footnotesize \obeylines \parindent=0pt

Daniela K\"{u}hn, Deryk Osthus 
School of Mathematics
University of Birmingham
Edgbaston
Birmingham
B15 2TT
UK
}
\begin{flushleft}
{\it{E-mail addresses}:
\tt{\{d.kuhn,d.osthus\}@bham.ac.uk}}
\end{flushleft}

\end{document}